\documentclass[12pt]{amsart}
\usepackage{a4wide}
\usepackage{graphicx}
\usepackage{epstopdf}
\usepackage{amssymb,amsmath, enumerate}
\usepackage{amsfonts}
\usepackage{psfrag}
\usepackage{graphicx,color}
\usepackage{amssymb}
\usepackage[pdfborder={0 0 0}]{hyperref}

\newtheorem{lem}{Lemma}[section]

\newtheorem{pro}{Proposition}[section]
\theoremstyle{definition}

\title{Infinite products of nonnegative $2\times2$ matrices by nonnegative vectors}
\author[A. Thomas]{Alain Thomas}
\address[Alain Thomas]{
LATP, 39, rue Joliot-Curie,\hfil\break
13453 Marseille, Cedex 13,
France}
\email{thomas@cmi.univ-mrs.fr}
\date{}  
\begin{document}
\baselineskip=18pt
\maketitle
\begin{abstract}Given a finite set $\{M_0,\dots,M_{d-1}\}$ of nonnegative $2\times 2$ matrices and a nonnegative column-vector $V$, we associate to each $(\omega_n)\in\{0,\dots,d-1\}^\mathbb N$ the sequence of the column-vectors $\displaystyle{M_{\omega_1}\dots M_{\omega_n}V\over\Vert M_{\omega_1}\dots M_{\omega_n}V\Vert}$. We give the necessary and sufficient condition on the matrices $M_k$ and the vector $V$ for this sequence to converge for all \hbox{$(\omega_n)\in\{0,\dots,d-1\}^\mathbb N$} such that $\forall n,\ M_{\omega_1}\dots M_{\omega_n}V\ne\begin{pmatrix}0\\0\end{pmatrix}$.
\end{abstract}
\par\par

\vskip20pt

\centerline{\sc Introduction}

\smallskip

Let  ${\mathcal M}=\{M_0,\dots,M_{d-1}\}$ be a finite set of nonnegative $2\times 2$ matrices and $V=\begin{pmatrix}v_1\\v_2\end{pmatrix}$ a nonnegative column-vector. We use the notation $Y_n=Y_n^\omega:=M_{\omega_1}\dots M_{\omega_n}$ and give the necessary and sufficient condition for the pointwise convergence of $\displaystyle{Y_nV\over\Vert Y_nV\Vert}$, \hbox{$(\omega_n)\in\{0,\dots,d-1\}^\mathbb N$} such that $Y_nV\ne\begin{pmatrix}0\\0\end{pmatrix}$ for any $n$, where $\Vert\cdot\Vert$ is the norm-sum. The idea of the proof is that, if the conditions are satisfied, either both columns of $Y_n$ tends to the same limit, or they tend to different limits with different orders of growth, so in case $V$ is positive the limit points of $\displaystyle{Y_nV\over\Vert Y_nV\Vert}$ only depend on the limit of the dominant column. This problem is obviously very different from the one of the convergence of $\displaystyle{Y_n\over\Vert Y_n\Vert}$, or the convegence of the $Y_n$ itselves, see the intoduction of~\cite{OTdXd} for some counterexamples and \cite[Proposition 1.2]{sto} for the infinite products of $2\times2$ stochastic matrices.

The conditions for the pointwise convergence of $\displaystyle{Y_nV\over\Vert Y_nV\Vert}$ also differ from the conditions for its uniform convergence, see \cite{OT2X2}. The uniform convergence can be used for the multifractal analysis of some continuous singular measures called Bernoulli convolutions (see \cite{PSS} for the Bernoulli convolutions and \cite{FO} for their multifractal analysis). We study such measures in \cite{OT2X2}, \cite{OT4X4} and \cite{OT7X7}. The Birkhoff's contraction coefficient \cite[Chapter 3]{S} that we use in \cite{OT4X4} and \cite{OTdXd} but not here, is really not of great help to solve the main difficulties. Moreover the theorem that gives the value of this coefficient is difficult to prove (see \cite[\S3.4]{S}) even in the case of $2\times2$ matrices. In~\cite{OT2X2} we use some other contraction coefficient quite more easy to compute (\cite[Proposition 1.3]{OT2X2}).

\section{Condition for the pointwise convergence of $\displaystyle{Y_nV\over\Vert Y_nV\Vert}$}

\begin{pro}{The sequence $\displaystyle{Y_nV\over\Vert Y_nV\Vert}$ converges for any $\omega\in\{0,\dots,d-1\}^\mathbb N$ such that $\forall n,\ Y_nV\ne\begin{pmatrix}0\\0\end{pmatrix}$, if and only if at least one of the following conditions holds:\hfill\break
(i) $V$ has positive entries and it is an eigenvector of any invertible matrix of the form $\begin{pmatrix}a&0\\0&d\end{pmatrix}$ or $\begin{pmatrix}0&b\\c&0\end{pmatrix}$ that belongs to $\mathcal M$.\hfill\break
(ii) Any invertible matrix $\begin{pmatrix}a&b\\c&d\end{pmatrix}\in{\mathcal M}$ satisfies $a>0$
and, if $b=c=0$, $a\ge d$.\hfill\break
(iii) Any invertible matrix $\begin{pmatrix}a&b\\c&d\end{pmatrix}\in{\mathcal M}$ satisfies $d>0$
and, if $b=c=0$, $d\ge a$.\hfill\break
(iv) $V$ has a null entry and all the invertible matrices $\begin{pmatrix}a&b\\c&d\end{pmatrix}\in{\mathcal M}$ satisfy $ad>0$.
}
\end{pro}
\begin{proof}Let $\omega\in\{0,\dots,d-1\}^\mathbb N$. If there exists N such that $\det M_{\omega_N}=0$, the column-vectors $Y_N V,\ Y_{N+1}V,\ \dots$ are collinear and $\displaystyle{Y_nV\over\Vert Y_nV\Vert}$ is constant for $n\ge N$.
So we look only at the $\omega\in\{0,\dots,d-1\}^\mathbb N$ such that $\forall n,\ \det
M_{\omega_n}\ne0$. In order to use only matrices with positive determinant we set $\Delta:=\begin{pmatrix}0&1\\1&0\end{pmatrix}$ and
$$
A_n=A_n^\omega:=\begin{cases}M_{\omega_n}\quad&\hbox{if }\det Y_{n-1}>0\ (\hbox{or }n=1)\hbox{ and }\det M_{\omega_n}>0\\
M_{\omega_n}\Delta\quad&\hbox{if }\det Y_{n-1}>0\ (\hbox{or }n=1)\hbox{ and }\det
M_{\omega_n}<0\\
\Delta M_{\omega_n}\quad&\hbox{if }\det Y_{n-1}<0\hbox{ and }\det M_{\omega_n}<0\\
\Delta M_{\omega_n}\Delta\quad&\hbox{if }\det Y_{n-1}<0\hbox{ and }\det
M_{\omega_n}>0.\end{cases}
$$
We set also
$$
\begin{pmatrix}a_n&b_n\\c_n&d_n\end{pmatrix}:= A_n\quad\hbox{and}\quad\begin{pmatrix}p_n&q_n\\
r_n&s_n\end{pmatrix}:=A_1\dots A_n
=\begin{cases}Y_n&\hbox{if }\det Y_n>0\\
Y_n\Delta&\hbox{if }\det Y_n<0.\end{cases}
$$
The matrices $A_n$ belong to the set
$$
{\mathcal M}^+:=\{M\;;\;\exists i,j,k,\ M=\Delta^iM_k\Delta^j\hbox{ and }\det M>0\}.
$$
Since $\det A_n>0$ we have $a_nd_np_ns_n\ne0$.
If $\{n;\ A_n\ \hbox{not\ diagonal}\}$ is infinite we index this set by an increasing sequence $n_1<n_2<\dots$.
We have $b_{n_1}\ne0$ or $c_{n_1}\ne0$; both cases are equivalent because, using the set of matrices ${\mathcal M}'=\Delta{\mathcal M}\Delta$ and defining similarly $Y'_n$ and $A'_n=\begin{pmatrix}a'_n&b'_n\\c'_n&d'_n\end{pmatrix}$ from this set,
we have
$Y'_n=\Delta Y_n\Delta$, $A'_n=\Delta A_n\Delta$ and
$b'_{n_1}=c_{n_1}$.

So we can suppose $b_{n_1}\ne0$; we deduce $q_n\ne0$ by induction on $n\ge n_1$. The sequences defined for any $n\ge n_1$ by
$$
u_n={r_n\over p_n},\quad
v_n={s_n\over q_n},\quad
w_n={q_n\over p_n},\quad x_n=\begin{cases}v_2/v_1&\hbox{if }\det
Y_n>0\\
v_1/v_2&\hbox{if not}\end{cases}\quad\hbox{and}\quad\lambda_n=(1+w_nx_n)^{-1}
$$
satisfy $0\le u_n<v_n<\infty$, $0<w_n<\infty$ and
$$
\begin{array}{l}\hbox{if the entries of }V\hbox{ are positive, }0<x_n<\infty\hbox{ and }0<\lambda_n<1\\
\hbox{if not, }x_n\in\{0,\infty\}\hbox{ and }\lambda_n\in\{0,1\}\hbox{ according to the sign of }\det Y_n.\end{array}
$$

Since we have assumed that $Y_nV\ne\begin{pmatrix}0\\0\end{pmatrix}$, the ratio $\displaystyle{(Y_nV)_2\over(Y_nV)_1}$ exists in $[0,\infty]$ and we have to prove that it has a finite or infinite limit when $n\to\infty$. If $A_n$ is not eventually diagonal we have for $n\ge n_1$
\begin{equation}\label{1}
{(Y_nV)_2\over(Y_nV)_1}=\lambda_nu_n+(1-\lambda_n)v_n\in I_n:=[u_n,v_n]\hbox{ and }I_n\supseteq I_{n+1}.
\end{equation}
An immediate consequence is the following lemma:

\begin{lem}\label{lem}Suppose $A_n$ is not eventually diagonal, then

(i) the sequences $(u_n)$ and $(v_n)$ converge in $\mathbb R$ and the sequence $\displaystyle\left({(Y_nV)_2\over(Y_nV)_1}\right)$ is bounded;

(ii) $\displaystyle\left({(Y_nV)_2\over(Y_nV)_1}\right)$ converges if $\displaystyle\lim_{n\to\infty}\vert I_n\vert=0$;

(iii) if $V$ has positive entries, $\displaystyle\left({(Y_nV)_2\over(Y_nV)_1}\right)$ converges if $\displaystyle w_n$ has limit $0$ or~$\infty$;

(iv) if $V$ has a null entry, the necessary and sufficient condition for the convergence of $\displaystyle\left({(Y_nV)_2\over(Y_nV)_1}\right)$ is that $\displaystyle\lim_{n\to\infty}\vert I_n\vert=0$ or the sign of $\det Y_n$ is eventually constant.
\end{lem}

We also define for $n>n_1$
$$
\alpha_n=\left(1+{c_n\over a_n}w_{n-1}\right)^{-1},\quad
\beta_n=\left(1+{b_n\over d_n}(w_{n-1})^{-1}\right)^{-1},\quad
\gamma_n=1-{c_n\over a_n}{b_n\over d_n}
$$
that belong to $]0,1]$ and satisfy
\begin{equation}\label{2}
\begin{array}{l}\displaystyle\vert I_n\vert=\alpha_n\beta_n\gamma_n\vert I_{n-1}\vert\\
\displaystyle w_n={d_n\over a_n}{\alpha_n\over\beta_n}w_{n-1}\end{array}
\end{equation}
so $\displaystyle\prod_{n>
n_1}\alpha_n\beta_n\gamma_n=\lim_{n\to\infty}{\vert I_n\vert\over\vert I_{n_1}\vert}$ is positive if and only if $\displaystyle\lim_{n\to\infty}\vert I_n\vert>0$. Using the equivalents of $\log\alpha_n$, $\log\beta_n$ and $\log\gamma_n$,
\begin{equation}\label{complicated}
\displaystyle\lim_{n\to\infty}\vert I_n\vert>0\Leftrightarrow\sum{c_n\over a_n}w_{n-1}<\infty,\quad
\sum{b_n\over d_n}(w_{n-1})^{-1}<\infty\quad\hbox{and}\quad
\sum{c_n\over a_n}{b_n\over d_n}<\infty.
\end{equation}
The set of indexes $\{n;\ A_n\ \hbox{not\ diagonal}\}$ is the union of
$$
L^\omega=\{n;\ c_n\ne0\}
\quad\hbox{and}\quad U^\omega=\{n;\ b_n\ne0\}.
$$
Moreover, since $A_n$ belongs to the finite set ${\mathcal M}^+$ there exists
$K>0$ such that
$$
L^\omega=\left\{n;\ {1\over K}\le{c_n\over a_n}\le K\right\}\quad\hbox{and}\quad U^\omega=\left\{n;\
{1\over K}\le{b_n\over d_n}\le K\right\}.
$$
We deduce a simpler formulation of (\ref{complicated}):
\begin{equation}\label{simple}
\displaystyle\lim_{n\to\infty}\vert I_n\vert>0\Leftrightarrow\sum_{n\in L^\omega}w_{n-1}<\infty,\quad\sum_{n\in U^\omega}(w_{n-1})^{-1}
<\infty\quad\hbox{and}\quad L^\omega\cap U^\omega\ \hbox{is finite}.
\end{equation}
In view of Lemma \ref{lem} we may suppose from now that $\displaystyle\lim_{n\to\infty}\vert I_n\vert>0$. Since $L^\omega\cap U^\omega=\{n\;;\;A_n\hbox{ positive}\}$ is finite, for $n$ large enough the matrix $A_n$ is lower triangular if $n\in L^\omega$, upper triangular if $n\in U^\omega$, diagonal if $n\not\in L^\omega\cup U^\omega$. When $A_n$ is diagonal the second relation in (\ref{2}) becomes $\displaystyle
w_n={d_n\over a_n}w_{n-1}$; consequently any integer $n$ in an interval $]n_i,n_{i+1}[$ with $i$ large enough satisfies
\begin{equation}\label{4}
{w_n\over w_{n_i}}=\prod_{n_i<j\le
n}{d_j\over a_j}.
\end{equation}
Moreover if $L^\omega$ is infinite, (\ref{simple}) implies that $w_{n-1}$ has limit to $0$ when $L^\omega\ni n\to\infty$, and $w_n$ also has limit $0$ because $\displaystyle w_n={d_nw_{n-1}\over a_n+c_nw_{n-1}}$ for any $n\in
L^\omega\setminus U^\omega$. We have a similar property if $U^\omega$ is infinite, so
\begin{equation}\label{5}
\begin{array}{ll}\hbox{if}\ L^\omega\ \hbox{is infinite},&
w_{n-1}\to0\ \hbox{and}\ w_n\to0\hbox{ for }L^\omega\ni n\to\infty;\\
\hbox{if}\ U^\omega\ \hbox{is infinite},&
w_{n-1}\to\infty\ \hbox{and}\ w_n\to\infty\hbox{ for }U^\omega\ni n\to\infty.\end{array}
\end{equation}

\smallskip

{\bf First case:} Suppose that (i) holds. Then the diagonal matrices
of $\mathcal M$
are collinear to the unit matrix. If at least one matrix of ${\mathcal M}$ has the form $M_k=\begin{pmatrix}0&b\\c&0\end{pmatrix}$ with $bc\ne0$, its nonnegative eigenvalue -- namely
$\begin{pmatrix}\sqrt b\\\sqrt c\end{pmatrix}$ -- is collinear to $V=\begin{pmatrix}v_1\\v_2\end{pmatrix}$ hence there exists $\lambda\in{\mathbb R}$ such that $M_k=\lambda\begin{pmatrix}0&v_1^2\\v_2^2&0\end{pmatrix}$.

Notice that if $A_n$ is diagonal from a rank $N$,
the matrix $M_{\omega_n}$ has the form $\begin{pmatrix}a&0\\0&d\end{pmatrix}$ or $\begin{pmatrix}0&b\\c&0\end{pmatrix}$ hence it has $V$ as eigenvector; consequently 
$\displaystyle\left({(Y_nV)_2\over(Y_nV)_1}\right)$ converges because it is $\displaystyle{(Y_NV)_2\over(Y_NV)_1}$ for any $n\ge N$.

Suppose now $A_n$ is non-diagonal for infinitely many $n$. We apply (\ref{4}) on each interval $]n_i,n_{i+1}[$ (if non empty), for $i$ large enough.
Among the integers $n\in]n_i,n_{i+1}[$ we consider the ones for which $\det M_{\omega_n}<0$.
For such $n$ the matrix $A_n$ is alternately
$M_{\omega_n}\Delta$ and $\Delta M_{\omega_n}$, hence alternately proportional to
$\begin{pmatrix}v_1^2&0\\0&v_2^2\end{pmatrix}$ and to $\begin{pmatrix}v_2^2&0\\0&v_1^2\end{pmatrix}$
and, according to (\ref{4}),
\begin{equation}\label{6}
n_i\le n<n_{i+1}\Rightarrow{w_n\over w_{n_i}}\in\Big\{{v_1^2\over
v_2^2},{v_2^2\over v_1^2},1\Big\}.
\end{equation}
In particular this relation holds for $n=n_{i+1}-1$. One deduce -- according to (\ref{5}) -- that there do not exist infinitely many $i$ such that
$n_i\in L^\omega$ and $n_{i+1}\in U^\omega$. Thus $n_i\in L^\omega$ for $i$ large enough (resp. $n_i\in U^\omega$ for $i$
large enough) and, according to (\ref{5}) and (\ref{6}), $\displaystyle\lim_{n\to\infty} w_n=0$ (resp. $\displaystyle\lim_{n\to\infty} w_n=\infty$). In view of Lemma \ref{lem}(iii),
$\displaystyle\left({(Y_nV)_2\over(Y_nV)_1}\right)$ converges.

\smallskip

{\bf Second case:} Suppose that (ii) holds
(if (iii) holds the proof is similar).

Suppose first the $M_{\omega_n}$ are diagonal from a rank $N$. From the hypothesis (ii) there exists $\delta_n,\delta'_n$ such that $M_{\omega_N}\dots M_{\omega_n}V=\left(\begin{array}{c}\delta_nv_1\\\delta'_nv_2\end{array}\right)$ and $\delta_n\ge\delta'_n$. Since the $M_{\omega_i}$ belong to a finite set we have $\displaystyle\lim_{n\to\infty}{\delta_n\over\delta'_n}=\infty$, or $\displaystyle{\delta_n\over\delta'_n}$ is eventually constant in case $M_{\omega_n}$ is eventually the unit matrix, or $\delta'_n=0\ne\delta_n$ for $n$ large enough. Denoting by $\left(\begin{array}{cc}p&q\\r&s\end{array}\right)$ the matrix $M_{\omega_1}\dots M_{\omega_{N-1}}$, we have $\displaystyle{(Y_nV)_2\over(Y_nV)_1}=\displaystyle{r\delta_nv_1+s\delta'_nv_2\over p\delta_nv_1+q\delta'_nv_2}$ hence $\displaystyle\left({(Y_nV)_2\over(Y_nV)_1}\right)$ converges in all the~cases.

Suppose now $M_{\omega_n}$ is non-diagonal for infinitely many $n$. There exists from (\ref{5}) an integer $\kappa$ such that
\begin{equation}\label{either}
i\ge\kappa\Rightarrow\left\{\begin{array}{ll}w_{n_i-1}<1\hbox{ and }w_{n_i}<1&\hbox{if }n_i\in L^\omega\\w_{n_i-1}>1\hbox{ and }w_{n_i}>1&\hbox{if }n_i\in U^\omega\end{array}\right.
\end{equation}
and such that the $A_n$ are diagonal for $n\in]n_i,n_{i+1}[$, $i\ge\kappa$. According to (ii), for such values of $n$ the matrix $M_{\omega_n}$ is diagonal and $A_n=M_{\omega_n}$ with $a_n\ge d_n$, or $A_n=\Delta M_{\omega_n}\Delta$ with $a_n\le d_n$.

If there exists $i\ge\kappa$ such that $n_i\in L^\omega$ and $n_{i+1}\in U^\omega$, $\det Y_{n_i}$ is necessarily negative: otherwise $A_n$ should be equal to $M_{\omega_n}$ for $n\in]n_i,n_{i+1}[$, $\displaystyle{d_n\over a_n}\le1$ and, by (\ref{4}), $w_{n_{i+1}-1}\le w_{n_i}<1$ in contradiction with (\ref{either}).

Now $M_{\omega_{n_{i+1}}}$ has positive determinant, otherwise it should have the form $\begin{pmatrix}a&b\\c&0\end{pmatrix}$ and $A_{n_{i+1}}=\Delta M_{\omega_{n_{i+1}}}=\begin{pmatrix}c&0\\a&b\end{pmatrix}$ in contradiction with $n_{i+1}\in U^\omega$.

We have again $\displaystyle{d_n\over a_n}\ge1$ for $n\in]n_{i+1},n_{i+2}[$ and consequently $w_{n_{i+2}-1}\ge w_{n_{i+1}}>1$; so, by induction, $n_j\in U^\omega$ and $\det Y_{n_j}<0$ for any $j\ge i+1$. From ({5}) $w_n$ lies between $w_{n_j}$ and $w_{n_{j+1}-1}$ for any $n\in]n_j,n_{j+1}[$ and $j$ large enough, and from (\ref{5}) its limit is infinite. Distinguishing the cases where $V$ has positive entries or $V$ has a null entry, $\displaystyle{(Y_nV)_2\over(Y_nV)_1}$ converges by Lemma \ref{lem}.

The conclusion is the same if there do not exist $i\ge\kappa$ such that $n_i\in L^\omega$ and $n_{i+1}\in U^\omega$, because in this case $n_i\in L^\omega$ for $i$ large enough, or $n_i\in U^\omega$ for any $i\ge\kappa$.

\smallskip

{\bf Third case:}
Suppose (iv) holds. As we have seen, from (\ref{simple}) $A_n$ is eventually triangular or diagonal, and $M_{\omega_n}$ also is because -- by (iv) -- $\mathcal M$ do not contain invertible matrices of the form $\begin{pmatrix}0&b\\c&d\end{pmatrix}$ or $\begin{pmatrix}a&b\\c&0\end{pmatrix}$. We deduce that the sign of $\det Y_n$ is eventually constant. If $A_n$ is not eventually diagonal the sequence $\displaystyle\left({(Y_nV)_2\over(Y_nV)_1}\right)$ converges by Lemma \ref{lem}(iv) and, if $A_n$~is, the sequence $\displaystyle\left({(Y_nV)_2\over(Y_nV)_1}\right)$ is eventually constant.

\smallskip

{\bf Fourth case:}
Suppose that the set $\mathcal M$ do not satisfy (i), (ii), (iii)
nor (iv), and that at~least one matrix of this set, let $M_k$, has the form $M_k=\begin{pmatrix}0&b\\c&0\end{pmatrix}$ with $bc\ne0$; let us prove that $\displaystyle\left({(Y_nV)_2\over(Y_nV)_1}\right)$ diverges.

$\_$ Suppose first there exists a matrix $M_k$ of this form that do not have $V$
as eigenvector; we chose as counterexample the constant sequence defined by $\omega_n=k$ for any $n$: $Y_{2n}$ is collinear to the unit matrix, hence $Y_{2n}V$ is collinear to $V$ and $Y_{2n+1}V$ to $M_kV$, so $\displaystyle\left({(Y_nV)_2\over(Y_nV)_1}\right)$ diverges.

$\_$ Suppose now that all the matrices of $\mathcal M$ of the form $\begin{pmatrix}0&b\\c&0\end{pmatrix}$ with $bc\ne0$ have $V$
as eigenvector that is, $V$ is collinear to $\begin{pmatrix}\sqrt b\\\sqrt c\end{pmatrix}$ for all such matrix. Since (i) do not hold, at least one matrix $M_h$ of $\mathcal M$ is diagonal with nonnull and distinct diagonal entries. In this case
$M_hM_k$ has the form $\begin{pmatrix}0&b\\c&0\end{pmatrix}$ but do not have $V$
as eigenvector. We recover the previous case; more precisely the counterexample is defined by
$\omega_{2n-1}=h$ and $\omega_{2n}=k$ for any $n\in\mathbb N$.

\smallskip

{\bf Fifth case:}
Suppose that $\mathcal M$ do not satisfy (i), (ii), (iii)
nor (iv),
and that no matrix of the form $\begin{pmatrix}0&b\\c&0\end{pmatrix}$ with $bc\ne0$ belongs to $\mathcal M$.
Since (i) do not hold, at least one matrix of this set has the form $M_k=\begin{pmatrix}a&0\\0&d\end{pmatrix}$ with $ad\ne0$ and $a\ne d$.
We suppose that $a>d$ and we use the negation of (ii) (in case $a<d$ we use similarly the negation of (iii)).
According to the negation of (ii) there exists in $\mathcal M$ at least one matrix of the form $M_h=
\begin{pmatrix}0&\beta\\\gamma&\delta\end{pmatrix}$ with $\beta\gamma\delta\ne0$,
or one of the form $M_\ell=
\begin{pmatrix}\alpha&0\\0&\delta\end{pmatrix}$ with $0<\alpha<\delta$.

$\_$ Consider first the case where $\mathcal M$ contains some matrices $M_k$ and $M_h$ as above.
Let $(n_i)_{i\in\mathbb N}$ be an increasing sequence of positive integers with $n_1=1$, and $\omega$ the sequence defined by
$\omega_n=h$ for $n\in\{n_1,n_2,\dots\}$ and $\omega_n=k$ otherwise.

For $i$ odd, $A_{n_i}$ is lower-triangular and $\forall n\in]n_i,n_{i+1}[,\ A_n=\begin{pmatrix}d&0\\0&a\end{pmatrix}$, $a_n=d$ and $d_n=a$. For $i$ even, $A_{n_i}$ is upper-triangular and $\forall n\in]n_i,n_{i+1}[,\ A_n=\begin{pmatrix}a&0\\0&d\end{pmatrix}$, $a_n=a$ and $d_n=d$.

Using (\ref{4}) for $n=n_{i+1}-1$ and choosing $n_{i+1}-n_i$ large enough one has $w_{n_{i+1}-1}\ge2^i$ if $i$ is odd, $w_{n_{i+1}-1}\le2^{-i}$ if $i$ is even, so the three conditions in (\ref{simple}) are satisfied and the interval $\displaystyle\cap I_n$ is not reduced to one point.
If the entries of $V$ are positive, the first relation in (\ref{1}) and the definition of~$\lambda_n$ imply that
$\displaystyle\liminf_{n\to\infty}\left({(Y_nV)_2\over(Y_nV)_1}\right)$ is the lower bound of this
interval and $\displaystyle\limsup_{n\to\infty}\left({(Y_nV)_2\over(Y_nV)_1}\right)$
its upper bound, so the sequence $\displaystyle\left({(Y_nV)_2\over(Y_nV)_1}\right)$ diverges. If $V$ has a null entry, the divergence of $\displaystyle\left({(Y_nV)_2\over(Y_nV)_1}\right)$ results from Lemma \ref{lem}(iv).

$\_$ In case $\mathcal M$ contains some matrices $M_k$ and $M_\ell$ as above, one defines $\omega$ from a sequence $i_1=1<i_2<i_3<\dots$
by setting, for $j\ge1$ and $i_j\le n<i_{j+1}$,
$$
\omega_n=\left\{\begin{array}{ll}k&\hbox{if }j\hbox{ even}\\\ell&\hbox{if }j\hbox{ odd}.\end{array}\right.
$$
The diagonal matrix $Y_n$ can be easily computed, and $\displaystyle\left({(Y_nV)_2\over(Y_nV)_1}\right)$ obviously diverges if one choose the
$i_{j+1}-i_j$ large enough.

If $V$ has a null entry, since (iv) do not hold $\mathcal M$ contains at least one matrix of the form $M_h=
\begin{pmatrix}0&\beta\\\gamma&\delta\end{pmatrix}$, $\beta\gamma\delta\ne0$ or $M_{h'}=\begin{pmatrix}\alpha&\beta\\\gamma&0\end{pmatrix}$, $\alpha\beta\gamma\ne0$, or $M_{h''}=\begin{pmatrix}0&\beta\\\gamma&0\end{pmatrix}$, $\beta\gamma\ne0$. We already know that $\displaystyle\left({(Y_nV)_2\over(Y_nV)_1}\right)$ diverges if $\mathcal M$ contains $M_k$ and $M_h$. Similarly it diverges if $\mathcal M$ contains $M_\ell$ and $M_{h'}$. If $\mathcal M$ contains $M_k$ and $M_{h''}$ the counterexample is given -- from a sequence $i_1=1<i_2<i_3<\dots$ -- by $\omega_{i_j}=h''$ and $\omega_n=k$ for $n\in]i_j,i_{j+1}[$, $j\in\mathbb N$:

$\displaystyle{(Y_nV)_2\over(Y_nV)_1}$ is alternately $0$ and $\infty$ because $Y_{i_j}$ has the form $\left(\begin{array}{cc}0&q\\r&0\end{array}\right)$ for $j$ odd and $\left(\begin{array}{cc}p&0\\0&s\end{array}\right)$ for $j$ even. 
\end{proof}

\end{document}